\documentclass[a4paper, twoside,10pt]{article}%
\usepackage{fancyhdr}
\usepackage{fnpos}
\usepackage[english]{babel}
\usepackage[T1]{fontenc}
\usepackage{graphicx}
\usepackage{makeidx}
\usepackage{fancybox}
\usepackage{framed}
\usepackage{fancyhdr}
\usepackage[margin=1in]{geometry}
\usepackage{graphicx}
\usepackage{titlesec}
\usepackage{amsmath}
\usepackage{amsfonts}
\usepackage{amssymb}
\usepackage{amsthm}
\usepackage{dsfont}
\usepackage{mathtools}
\usepackage{verbatim}
\usepackage{color}
\usepackage{listings}
\usepackage{graphicx}
\usepackage{amsmath}
\usepackage{amsmath,amssymb,color,inputenc,euscript,graphicx,psfrag}
\usepackage{amsmath, amsthm, amscd, amsfonts, amssymb, graphicx, color,
mathrsfs}
\usepackage[english]{babel}
\usepackage[bookmarksnumbered, colorlinks, plainpages]{hyperref}%
\providecommand{\U}[1]{\protect\rule{.1in}{.1in}}

\hypersetup{colorlinks=true,linkcolor=blue, anchorcolor=blue, citecolor=blue, urlcolor=red, filecolor=magenta, pdftoolbar=true}
\newtheorem{theorem}{Theorem}[section]
\theoremstyle{plain}

\newtheorem{lemma}{Lemma}[section]

\newtheorem{remark}{Remark}[section]

\numberwithin{equation}{section}

\newcommand{\R}{\mathbb{R}}

 {\markboth{\sl  J. Benameur}{\sl Decay rates for mild solutions of QGE ... }
\begin{document}

\author{ Jamel Benameur$^{a}$}
\title{Decay rates for mild solutions of QGE with critical fractional dissipation in $L^2(\R^2)$}
\date{}
\maketitle

\begin{abstract}
In \cite{MRSC1} the authors proved some asymptotic results for the global solution of critical Quasi-geostrophic equation with a condition on the decay of $\widehat{\theta_0}$ near at zero. In this paper, we prove that this condition is not necessary. Fourier analysis and standard techniques are used.
\end{abstract}

\begin{center}
$^{a}$ University of Gabes, Department of Mathematics, Faculty of Sciences,
Tunisia.\newline\textbf{ e-mail:} jamelbenameur@gmail.com
\end{center}

\section{Introduction}
In this paper, we consider the Quasi-Geostrophic Equations with critical fractional dissipation defined on $\R^2$ by:
$$(1)\left\{\begin{array}{l}
\partial_t\theta+|D|\theta+u_\theta\nabla \theta=0\,in\,\R^+\times\R^2\\
\theta(0)=\theta^0\,in\,\R^2,
\end{array}\right.$$
where $\theta(x, t)$ denotes the potential temperature, and $u_\theta(x, t)$ stands for the fluid
velocity field. Here, we define
$$u_\theta=(\partial_2(-\Delta)^{1/2}\theta,-\partial_1(-\Delta)^{1/2}\theta)$$
and, as a result, one infers
$${\rm div}\,u_\theta=0\,{\rm and}\,|\widehat{u_\theta}|=|\widehat{\theta}|.$$
Furthermore, the fractional Laplacian $(-\Delta)^{1/2}=|D|$ is defined by
$$(-\Delta)^{1/2}f=\mathcal F^{-1}(|\xi|\widehat{f}(\xi))$$ for
all $f\in S'(\R^2)$ such that $\widehat{f}$ is given by a local integrable function. If $\theta_0$ is more regular, we have the Duhamel formula
\begin{equation}\label{deq1}\theta(t,x)=e^{-t|D|}\theta_0(x)-\int_0^te^{-(t-z)|D|}{\rm div}\,(\theta u_\theta)(z,x)dz.\end{equation}
Now we recall the results of \cite{MRSC1}:
\begin{theorem}\label{th1}(\cite{MRSC1}) Let $a > 0$, $\sigma > 1$, and $0\leq s < 1$. Assume that $\theta_0\in\dot H^s_{a,\sigma}(\R^2)$. There exists a
positive constant $C_{a,s,\sigma}$ such that if $\|\theta^0\|_{\dot H^s_{a,\sigma}}< C_{a,s,\sigma}$; then, one obtains a unique solution
for the Quasi-geostrophic equation (1), $\widetilde{L^\infty}(\R^+,\dot H^s_{a,\sigma}(\R^2))\cap L^2(\R^+,\dot H^{s+\frac{1}{2}}_{a,\sigma}(\R^2))$.
Moreover, the analyticity of this solution is related to the following inequality:
\begin{equation}\label{R01}\|\theta\|_{\widetilde{L^\infty}(\R^+,\dot H^s_{a,\sigma}(\R^2))}+\|\theta\|_{L^2(\R^+,\dot H^{s+\frac{1}{2}}_{a,\sigma}(\R^2))}\leq6\|\theta_0\|_{\dot H^s_{a,\sigma}(\R^2)}.\end{equation}
Moreover, we have
\begin{equation}\label{R02}\|e^{t^{1/2}|D|^{1/2}}\theta\|_{\widetilde{L^\infty}(\R^+,\dot H^s_{a,\sigma}(\R^2))}+\|\theta\|_{L^2(\R^+,\dot H^{s+\frac{1}{2}}_{a,\sigma}(\R^2))}\leq6e\|\theta_0\|_{\dot H^s_{a,\sigma}(\R^2)}.\end{equation}
\end{theorem}
\begin{theorem}\label{th2}(\cite{MRSC1})
Let $a>0$, $\sigma>1$, $0\leq s<1$, $\theta_0\in\dot H^s_{a,\sigma}(\R^2)\cap L^2(\R^2)$, and suppose that $r^*(\theta_0)$ exists and $r^*(\theta_0)=(-1,+\infty)$. Thus, there is a positive constant $C_{a,\sigma,s}'$ such that if $\|\theta^0\|_{\dot H^s_{a,\sigma}}<C_{a,\sigma,s}'$, then, the solution $\theta$ obtained in Theorem \ref{th1} satisfies the following statements:
\begin{enumerate}
\item[{\bf i)}] If $s=0$, one concludes that
$$\limsup_{t\rightarrow\infty} t^\kappa \|\theta(t)\|_{\dot H^\kappa_{a,\sigma}}=0,\,\forall \kappa\geq0,$$
\item[{\bf ii)}] If $s\in(0,1)$, we deduce that
$$\limsup_{t\rightarrow\infty} t^\kappa \|\theta(t)\|_{\dot H^\kappa_{a,\sigma}}=0,\,\forall \kappa>0,$$
\end{enumerate}
where $r^*(\theta_0)$ is the decay character for $\theta^0$ (for more details of the definition of decay character, see \cite{NS1} and references therein).
\end{theorem}

{\bf The conditions $r^*(\theta)$ exists and  $r^*(\theta)=(-1,+\infty)$ are not necessary to prove Theorem\ref{th2}-i)-ii).} Precisely, we have the following result.
\begin{theorem}\label{th3} Let $a>0$, $\sigma>1$, $0\leq s<1$, $\theta_0\in\dot H^s_{a,\sigma}(\R^2)\cap L^2(\R^2)$. Thus, there is a positive constant $C_{a,\sigma,s}'$ such that if $\|\theta^0\|_{\dot H^s_{a,\sigma}}<C_{a,\sigma,s}'$, then, the solution $\theta$ obtained in Theorem \ref{th1} satisfies the following statements:
\begin{equation}\label{R1}\limsup_{t\rightarrow\infty} \|\theta(t)\|_{L^2}=0\end{equation}
\begin{equation}\label{R2}\limsup_{t\rightarrow\infty} t^\kappa \|\theta(t)\|_{\dot H^\kappa_{a,\sigma}}=0,\,\forall \kappa\geq0.\end{equation}
\end{theorem}
\begin{remark}
\begin{enumerate}
\item If $\theta_0\in\dot H^s_{a,\sigma}(\R^2)\cap L^2(\R^2)$, then the solution $\theta$ of the system $(1)$ given by Theorem \ref{th1} satisfies the $L^2-$ energy estimate
\begin{equation}\label{eeq}\|\theta(t)\|_{L^2}^2+2\int_0^t\||D|^{1/2}\theta(z)\|_{L^2}^2=\|\theta_0\|_{L^2}^2,\,\forall t\geq0.\end{equation}
\item Inequalities (\ref{R02}) and (\ref{eeq}) imply that $\theta(t)\in H^\kappa_{a,\sigma}(\R^2)$ for all $t>0$ and $\kappa\geq0$. Indeed: For $t>0$. We distinguish two cases\\
    First case $\kappa\geq s$. By using inequality (\ref{R02}), we get
    $$\begin{array}{lcl}
    \displaystyle\int_{\R^2}|\xi|^{2\kappa}e^{2a|\xi|^{1/\sigma}}|\widehat{\theta}(t,\xi)|^2d\xi&=&
    \displaystyle\int_{\R^2}|\xi|^{2(\kappa-s)}
    e^{-2t^{1/2}|\xi|^{1/2}}|\xi|^{2s)}e^{2t^{1/2}|\xi|^{1/2}}e^{2a|\xi|^{1/\sigma}}|\widehat{\theta}(t,\xi)|^2d\xi\\
    &\leq& C(t,\kappa,s)
    \displaystyle\int_{\R^2}|\xi|^{2s)}e^{2t^{1/2}|\xi|^{1/2}}e^{2a|\xi|^{1/\sigma}}|\widehat{\theta}(t,\xi)|^2d\xi\\
    &\leq& C(t,\kappa,s)64e^2\|\theta_0\|_{\dot H^s_{a,\sigma}}^2,
    \end{array}$$
    where $C(t,\kappa,s)=\sup_{z\geq0}z^{2(\kappa-s)}
    e^{-2t^{1/2}z^{1/2}}<\infty$.\\
 Second case $0\leq\kappa< s$. We have
    $$\begin{array}{lcl}
    \displaystyle\int_{\R^2}|\xi|^{2\kappa}e^{2a|\xi|^{1/\sigma}}|\widehat{\theta}(t,\xi)|^2d\xi&=&
    \displaystyle\int_{|\xi|\leq1}|\xi|^{2\kappa}e^{2a|\xi|^{1/\sigma}}|\widehat{\theta}(t,\xi)|^2d\xi+
    \int_{|\xi|>1}|\xi|^{2\kappa}e^{2a|\xi|^{1/\sigma}}|\widehat{\theta}(t,\xi)|^2d\xi\\
    &\leq&
    \displaystyle e^{2a}\int_{|\xi|\leq1}|\widehat{\theta}(t,\xi)|^2d\xi+
    \int_{|\xi|>1}|\xi|^{2\kappa}e^{2a|\xi|^{1/\sigma}}|\widehat{\theta}(t,\xi)|^2d\xi\\
&\leq&\displaystyle e^{2a}\|\widehat{\theta}(t)\|_{L^2}^2+
    \int_{|\xi|>1}|\xi|^{2s}e^{2a|\xi|^{1/\sigma}}|\widehat{\theta}(t,\xi)|^2d\xi\\
&\leq&\displaystyle e^{2a}\|\widehat{\theta}(t)\|_{L^2}^2+
    \|\theta(t)\|_{\dot H^s_{a,\sigma}}^2\\
    &\leq&e^{2a}\|\widehat{\theta_0}\|_{L^2}^2+36\|\theta_0\|_{\dot H^s_{a,\sigma}}^2.
    \end{array}$$
\item The smoothing effect of the system (1) is given by the strong result (\ref{R02}), which gives the asymptotic results of Theorem \ref{th2}.
\item To prove Theorem \ref{th3}, it suffices to prove (\ref{R1}) and if $0<s<1$ the following limit
\begin{equation}\label{R3}\limsup_{t\rightarrow\infty} \|\theta(t)\|_{\dot H^0_{a,\sigma}}=0.\end{equation}
The study of $$\limsup_{t\rightarrow\infty} t^\kappa \|\theta(t)\|_{\dot H^\kappa_{a,\sigma}},\,\forall \kappa>0$$ was carried out in the reference \cite{MRSC1}.
\end{enumerate}
\end{remark}
To prove this theorem, we need the following technical lemmas.
\begin{lemma}\label{lem1}(\cite{BCD1}) Let $s_1,\ s_2$ be two real numbers.
\begin{enumerate}
\item If $s_1<1$\; and\; $s_1+s_2>0$, there exists a constant
$C_1=C_1(s_1,s_2)$, such that: if $f,g\in
\dot{H}^{s_1}(\mathbb{R}^2)\cap \dot{H}^{s_2}(\mathbb{R}^2)$, then
$f.g \in \dot{H}^{s_1+s_2-1}(\mathbb{R}^2)$ and
$$\|fg\|_{\dot{H}^{s_1+s_2-1}}\leq C_1 (\|f\|_{\dot{H}^{s_1}}\|g\|_{\dot{H}^{s_2}}+\|f\|_{\dot{H}^{s_2}}\|g\|_{\dot{H}^{s_1}}).$$
\item If $s_1,s_2<1$\; and\; $s_1+s_2>0$ there exists a constant
$C_2=C_2(s_1,s_2)$ such that: if $f \in
\dot{H}^{s_1}(\mathbb{R}^3)$\; and\;
$g\in\dot{H}^{s_2}(\mathbb{R}^2)$, then  $f.g \in
\dot{H}^{s_1+s_2-\frac{3}{2}}(\mathbb{R}^2)$ and
$$\|fg\|_{\dot{H}^{s_1+s_2-1}}\leq C_2 \|f\|_{\dot{H}^{s_1}}\|g\|_{\dot{H}^{s_2}}.$$
\end{enumerate}
\end{lemma}
\begin{lemma}\label{lem2}\cite{BK1} Let $h:[0,T]\rightarrow[0,\infty]$ be measurable function and $\sigma>0$, then
$$\Big(\int_0^Te^{-\sigma(T-z)}h(z)dz\Big)^2\leq 2\sigma^{-1}\int_0^Te^{-\sigma(T-z)}h(z)^2dz.$$
\end{lemma}
\section{Proof of Theorem\ref{th3}}
\begin{enumerate}
\item[$\bullet$] Proof of (\ref{R1}): The idea of this proof is due to \cite{J1} and \cite{BK1} which is based to an interpolation between positive and negative homogeneous Sobolev spaces.
We have
$$\begin{array}{lcl}
\displaystyle\int_0^\infty\|\int_0^te^{-(t-z)|D|}u_\theta\nabla \theta\|_{\dot H^{-1/2}}^2dt&=&\displaystyle\int_0^\infty\|\int_0^te^{-(t-z)|D|}{\rm div\,}(\theta u_\theta)\|_{\dot H^{-1/2}}^2dt\\
&\leq&\displaystyle\int_0^\infty\int_{\R^2}|\xi|^{-1}\big(\int_0^te^{-(t-z)|\xi|}|\xi|.|\mathcal F(\theta u_\theta)(z,\xi)dz|\Big)^2d\xi dt\\
&\leq&\displaystyle\int_0^\infty\int_{\R^2}|\xi|\Big(\int_0^te^{-(t-z)|\xi|}.|\mathcal F(\theta u_\theta)(z,\xi)dz|\Big)^2d\xi dt.
\end{array}$$
By using Lemma \ref{lem2}, we get
$$\begin{array}{lcl}
\displaystyle\int_0^\infty\|\int_0^te^{-(t-z)|D|}u_\theta\nabla \theta\|_{\dot H^{-1}}^2dt
&\leq&\displaystyle\int_0^\infty\int_{\R^2}\int_0^te^{-(t-z)|\xi|}.|\mathcal F(\theta u_\theta)(z,\xi)|^2dz d\xi dt\\
&\leq&\displaystyle\int_{\R^2}\Big(\int_0^\infty\int_0^te^{-(t-z)|\xi|}.|\mathcal F(\theta u_\theta)(z,\xi)|^2dzdt\Big) d\xi\\
&\leq&\displaystyle\int_{\R^2}\Big(\int_0^\infty(\int_z^\infty e^{-(t-z)|\xi|}dt).|\mathcal F(\theta u_\theta)(z,\xi)|^2dz\Big) d\xi\\
&\leq&\displaystyle\int_{\R^2}\Big(\int_0^\infty(|\xi|^{-1}).|\mathcal F(\theta u_\theta)(z,\xi)|^2dz\Big) d\xi\\
&\leq&\displaystyle\int_0^\infty(\int_{\R^2}|\xi|^{-1}.|\mathcal F(\theta u_\theta)(z,\xi)|^2d\xi dz\\
&\leq&\displaystyle\int_0^\infty\|\theta u_\theta\|_{\dot H^{-1/2}}^2dt
\end{array}$$
Using the product law in homogeneous Sobolev spaces(see Lemma \ref{lem1}) with $s_1=0$ and $s_2=1/2$, we get
$$\int_0^\infty\|\int_0^te^{-(t-z)|D|}u_\theta\nabla \theta\|_{\dot H^{-1}}^2dt
\leq\displaystyle C\int_0^\infty \|\theta\|_{L^2}^2\|\theta\|_{\dot H^{1/2}}^2dt.$$
The energy estimate (\ref{eeq}) implies that
$$\int_0^\infty\|\int_0^te^{-(t-z)|D|}u_\theta\nabla \theta\|_{\dot H^{-1}}^2dt
\leq\displaystyle C'\|\theta^0\|_{L^2}^4,$$
which gives that $$w:=\Big(t\longmapsto \int_0^te^{-(t-z)|D|}u_\theta\nabla \theta\Big)\in L^2(\R^+,\dot H^{-1/2}(\R^2)).$$
But $\gamma:=\Big(t\longmapsto e^{-t|D|}\theta^0\Big)\in L^2(\R^+,\dot H^{1/2}(\R^2))$ and $\theta\in L^2(\R^+,\dot H^{1/2}(\R^2))$, then
$$w \in L^2(\R^+,\dot H^{1/2}(\R^2)).$$
By interpolation, we get
$$\|w\|_{L^2}^2=c_0\|w\|_{\dot H^0}^2\leq c_0\|w\|_{\dot H^{-1/2}}\|w\|_{\dot H^{1/2}},$$
which implies, that $w\in L^2(\R^+,L^2(\R^2))$.\\
By Dominated Convergence Theorem, we obtain
$$\lim_{t\rightarrow\infty}\|e^{-t|D|}\theta^0\|_{L^2}=0.$$Then for a fixed real number $\varepsilon>0$, there is a time $t_0\geq0$ such that
\begin{equation}\label{eq11}\|e^{-t|D|}\theta^0\|_{L^2}<\frac{\varepsilon}{2},\,\forall t\geq t_0.\end{equation}
As $w\in L^2(\R^+,L^2(\R^2))$, then $w\in L^2([t_0,\infty),L^2(\R^2))$ and
$$\lambda_1(A_\varepsilon)<\infty,$$
where
$$A_\varepsilon:=\{t\geq t_0:\,\|w(t)\|_{L^2}\geq \varepsilon/2\}.$$
Therefore, there is $t_1\geq t_0$ such that
\begin{equation}\label{eq12}\|w(t_1)\|_{L^2}< \varepsilon/2.\end{equation}
By  (\ref{eq11}) and (\ref{eq12}) and the Duhamel formula (\ref{deq1}), we get
$$\|\theta(t_1)\|_{L^2}\leq\|\gamma(t_1)\|_{L^2}+\|w(t_1)\|_{L^2}< \varepsilon.$$
By the uniqueness of solution of the following system(given by Theorem \ref{th1})
$$(2)\left\{\begin{array}{l}
\partial_t\phi+|D|\phi+u_\phi\nabla \phi=0\\
\phi(0)=\theta(t_1)
\end{array}\right.$$
we get $\phi=\theta(t_1+t)$ is the solution of (2). Then by $L^2-$energy estimate (\ref{eeq}), we get for all $t\geq 0$
$$\|\phi(t)\|_{L^2}\leq \|\phi(0)\|_{L^2},\,\forall t\geq 0$$
and
$$\|\theta(t_1+t)\|_{L^2}\leq \|\theta(t_1)\|_{L^2}<\varepsilon,\,\forall t\geq 0.$$
which compete the proof of (\ref{R1}).
\item[$\bullet$] Proof of (\ref{R3}) : It suffices to suppose $0<s<1$, the other cases are treated in Theoerm \ref{th2}. For $t\geq 1$, we have
$$\begin{array}{lcl}
\|\theta(t)\|_{\dot H^0_{a,\sigma}}^2&=&\displaystyle\int_{\R^2}e^{2a|\xi|^{1/\sigma}}|\theta(t,\xi)|^2d\xi\\
&\leq&\displaystyle\int_{|\xi|<1}e^{2a|\xi|^{1/\sigma}}|\theta(t,\xi)|^2d\xi+\int_{|\xi|>1}e^{2a|\xi|^{1/\sigma}}|\theta(t,\xi)|^2d\xi\\
&\leq&\displaystyle e^{2a}\int_{|\xi|<1}|\theta(t,\xi)|^2d\xi+\int_{|\xi|>1}|\xi|^{2}e^{2a|\xi|^{1/\sigma}}|\theta(t,\xi)|^2d\xi\\
&\leq&\displaystyle e^{2a}\|\widehat{\theta}(t)\|_{L^2}^2+\|\theta(t)\|_{\dot H^1_{a,\sigma}}^2.
\end{array}$$
By the result proved in Theorem \ref{th2}-{\bf ii)} for $\kappa=1$, we get
$$\lim_{t\rightarrow\infty}\|\theta(t)\|_{\dot H^1_{a,\sigma}}^2=0.$$
Combining this result with (\ref{R1}) we get
$$\lim_{t\rightarrow\infty}\|\theta(t)\|_{\dot H^0_{a,\sigma}}^2=0.$$
\end{enumerate}
Then, the proof of Theorem \ref{th3} is completed.

\end{document}